\documentclass[a4paper,reqno,10pt]{amsart}
\usepackage{amsmath}
\usepackage{amsthm}
\usepackage{amssymb}
\usepackage{amscd} %diagrams
\usepackage{stmaryrd}

\usepackage{bbm} %bold numbers
\usepackage{mathrsfs} %nice calligraphy
\usepackage{verbatim} %DO NOT comment this line
\newtheoremstyle{exercise} %for books or class notes
  {3pt} %space above
  {3pt} %space below
  {\scriptsize\rmfamily} %body font
  {
\parindent} %indent amount(empty=no indent,\parindent=para indent)
  {\rmfamily\scshape} %thm head font
  {.} %punctuation after thm head
  {.5em} %space after thm head: " " = normal interword space;
  {} %thm head spec (can be left empty, meaning `normal')

\newtheoremstyle{newplain}
  {5pt}
  {5pt}
  {\itshape}
  {}
  {\rmfamily\scshape}
  {. ---}
  {.5em}
  {}

\newtheoremstyle{newremark}
  {5pt}
  {5pt}
  {\rmfamily}
  {}
  {\rmfamily\scshape}
  {. ---}
  {.5em}
  {}

\theoremstyle{newplain} 
\newtheorem{theorem}{Theorem}[section] 
\newtheorem{lemma}[theorem]{Lemma} 
 
\newtheorem{corollary}[theorem]{Corollary} 
\newtheorem{definition}[theorem]{Definition} 
 
\newtheorem{remark}[theorem]{Remark} 

\newcommand{\Nbb}{\mathbb{N}} 
\newcommand{\Tbb}{\mathbb{T}}

\newcommand\loc{\operatorname{\loc}}

\newcommand{\one}{\hbox{1\kern -2.4pt 1}}

\newcommand\s{\sigma}

\newcommand{\N}{\Nbb} 
 
\def\tpsi{\tilde{\psi}}

\def\xchi{\smash{\raise 0.4ex\hbox{$\chi$}}} 
\def\<{\langle} 
\def\>{\rangle}

\def\Neal{\Bbb {I\kern-2pt N}} 
\def\neal{\Bbb {I\kern-2pt N}}

\def\pf{\begin{proof}}
\def\thst{ }

\begin{document} 
\title[Ergodic transformations generate Carath\'eodory's definition]{Outer measure preserving ergodic transformations generate the Carath\'eodory definition of measurable sets} 
{\small \author{Amos N. Koeller }}

\thanks{This research was supported by a University Postgraduate Scholarship (UPS) at the University of Wollongong, 
Australia. 
}

\begin{abstract} \noindent
It is known that there are specific examples of ergodic transformations on measure spaces for which the calculation of 
the outer measure of transformation invariant sets leads to a condition closely resembling Carath\'eodory's condition 
for sets to be measurable. It is then natural to ask what functions are capable of `generating', that is leading to, the Carath\'eodory definition in the same way. The present work answers this question by showing that the property 
of generating Carath\'eodory's definition holds for the general class of outer measure preserving ergodic 
transformations on measure spaces. We further show that the previously found specific examples of functions generating Carath\'eodory's definition fall into this family of transformations.
\end{abstract}
\maketitle 
\section {Introduction} 
\setcounter{equation}{0} 
In the theory of measure and integration the definition given by Carath\'eodory \cite{Ca} in 1914 is a crucial point in the relationship between an outer measure defined on the set of all subsets of some given set, and a measure defined on 
some $\s$-algebra of subsets for which additional convenient properties with respect to the measure hold 
(~\cite[pp.101-103]{Ba} or ~\cite[pp.44-46]{Ha}).
Carath\'eodory's definition is that a subset $A$ of some set $X$ on which an outer measure $\mu_{\ast}$ is defined is called 
{\it measurable} if for all subsets $B$ of $X$ 
\begin{equation}
\mu_{\ast}(B) =\mu_{\ast}(B \cap A)+\mu_{\ast}(B\cap A^c)
\label{ecaradef}
\end{equation}
(where from here on $A^c$ denotes the complement of the set $A$). 

The apparent lack of sufficient motivation for this definition has been commented upon by several mathematicians (see, for example, the discussion in \cite{Ni1}). In \cite{Ni1}, however, it was shown that the definition of Carath\'eodory can be seen to naturally arise by calculating the outer measure of the invariant sets of an irrational rotation on the unit circle. Specifically, let $\Tbb$ denote the set of complex numbers of modulus 1, and note that an irrational rotation on $\Tbb$ is an ergodic, but not weakly mixing, transformation of the form $z\mapsto \rho z$. Letting $\mu$ denote the usual Borel measure on $\Tbb$ and $\mu_{\ast}$ be the associated outer measure, it was proven in \cite{Ni1} that if $\rho$ is an irrational rotation, $A$ is a $\rho$-invariant subset of $\Tbb$, and there is a $\theta<2$ such that for all arcs $J$ of $\Tbb$, 
\begin{equation}
\mu_{\ast}(A \cap J)+\mu_{\ast}(A^c \cap J) \leq \theta \mu(J)
\label{reasonable}
\end{equation}
then either $\mu_{\ast}(A)=0$ or $\mu_{\ast}(A^c)=0$. The simplest way for ($\ref{reasonable}$) to hold is for the set $A$ to satisfy the Carath\'eodory condition, and in this sense, as motivated more fully in \cite{Ni1}, the problem has motivated Carath\'eodory's definition of measurable sets.

Having shown that Carath\'eodor's definition can be generated by some transformations, we ask whether this is a coincidence, or whether this is a property of some large class of transformations for which such a connection makes sense. Candidates for such a class should be sufficiently measure preserving and should take any test set of positive measure to any other. That is, the obvious class of transformations would be the class of measure preserving transformations. In \cite{Ni1} it was conjectured that the class may be the set of measure preserving ergodic transformations that are not strongly mixing. The later clarification due to the additional conjecture that the transformation $z\mapsto z^2$ on $\Tbb$ does not generate Carath\'eodory's definition in this way.

In \cite{Ko} it was shown that a weakly, but not strongly, mixing ergodic transformation due to Kakutani \cite{Ka2} generates Carath\'eodory's definition in the same sense as in \cite{Ni1}, which appeared to support the conjecture in \cite{Ni1}.

In this paper, however, we show that the family of transformations generating Carath\'eodory's definition in the sense described above is the more general class of all outer measure preserving ergodic transformations. We go on to show that each of the previously individually treated examples fall into this class, including the transformation $z \mapsto z^2$ on $\Tbb$. We therefore answer the conjectures raised in \cite{Ni1} and find the desired large, heuristically expected class that can be considered to generate Carath\'eodory's definition.

\section{Definitions and the analytic principles of generating Carath\'eodory's definition}
\setcounter{equation}{0}
The idea in both the irrational rotations and Kakutani's transformation cases was to take basic sets of measure estimation (in 
the former case arcs and in the latter case dyadic intervals, see \cite{Ko2}), say $J_1$ and $J_2$, 
of the same measure to find areas of density of the invariant set $A$ and its complement $A^c$. It was shown that 
$J_1$ and $J_2$ could be found such that 
\begin{equation}
\mu_{\ast}(A\cap J_1)>(1-\varepsilon)\mu(J_1) \hbox{ and }\mu_{\ast}(A^c\cap J_2)>(1-\varepsilon)\mu(J_2)
\label{thicksets}
\end{equation} 
for $\varepsilon$ as small as desired. It was then shown 
(using the properties of the transformation in question) that either $J_1$ or $J_2$ could, under negative iterations of 
the transformations, be `rotated' onto the other (in the irrational rotations case not exactly onto the other, but 
as nearly as required). Using the measure preserving properties of the transformations and ($\ref{thicksets}$), 
a contradiction to ($\ref{reasonable}$) would then arise.

We now extend these ideas to a general measure space $(X,\mathcal{B},\mu)$ equipped with an outer measure 
preserving transformation $T$ (together forming the dynamical system $(X,\mathcal{B},\mu,T)$). 

In order to give more general results, we need first to establish what it means, in general, to generate Carath\'eodory's definition. We consider firstly the role of the arcs or intervals in the previous works, which was to be a simple family of measurable sets from which, for any subset of $A\subset X$, an element could be found which locally represented the measure of $A$ well. We call such a family, which is somewhat like a topological subbasis, a measure basis. A Definition of a measure basis is given below in Definition $\ref{scos}$ after the concept of representing the measure of a set well is defined as a set of density in Definition $\ref{setsofdensity}$.

\begin{definition}\label{setsofdensity}\thst
If $(X,\mathcal{B},\mu)$ is a measure space, $B \subset X$ has positive outer measure and $\varepsilon>0$, then 
$J \in \mathcal{B}$ is called a set of density to within $\varepsilon$ for $B$ if 
$\mu_{\ast}(B\cap J)>(1-\varepsilon)\mu(J).$
\end{definition}\noindent

\begin{definition}\label{scos}\thst
Let $(X, \mathcal{B},\mu)$ be a measure space. A collection $\mathcal{J} \subset \mathcal{B}$ will be called a {\it measure basis} if for each pair of sets of positive measure $A_1, A_2 \subset X$ and each $\varepsilon>0$ there 
exist sets $J_1, J_2 \in \mathcal{J}$ of equal measure such that $J_i$ is a set of density to within 
$\varepsilon$ for $A_i$ for each $i\in \{1,2\}$.
\end{definition}\noindent 

\begin{remark}
\begin{enumerate}
\item We note that in any non-atomic measure space we will in fact be able to find a set in $\mathcal{B}$, say $J$, for 
which $\mu_{\ast}(B\cap J)=\mu(J)$. We do not dwell on this fact though, as we wish to allow for the selection of as small a measure basis as possible.

\item Note that every dynamical system will have at least one measure basis as $\mathcal{B}$ 
itself is a measure basis. Note also that the set of arcs in $\Tbb$ and the set of dyadic intervals in $[0,1)$ as considered in \cite{Ni1} and \cite{Ko} respectively are measure bases.
\end{enumerate}
\end{remark}

Having defined the basic objects to be considered in this work, we can now define how a transformation is considered to generate the Carath\'eodory definition analogously to the interpretations of \cite{Ko} and \cite{Ni1}.

\begin{definition}\label{generate}\thst
Let $T$ be a $\mu$-outer measure preserving transformation on a measure space $(X,\mathcal{B},\mu)$. That is, a transformation $T:X\rightarrow X$ such that $\mu_{\ast}(T^{-1}(A))=\mu_{\ast}(A)$ for all $A\subset X$. Furthermore, 
let $\mathcal{J}$ be a measure basis for $(X,\mathcal{B},\mu)$. Then, if whenever $A$ is a $T$-invariant subset of $X$ for which there is a 
$\theta=\theta(A)<2$ such that for all $J \in \mathcal{J}$
$$\mu_{\ast}(A\cap J) + \mu_{\ast}(A^c \cap J) \leq \theta\mu(J),$$
either $A$ or $A^c$ is a set of outer measure zero, $T$ is said to generate the Carath\'eodory definition of 
measurable sets relative to $\mathcal{J}$ or simply to generate Carath\'eodory's definition.
\end{definition}\noindent

\section{Transformations generating Carath\'eodory's definition}
\setcounter{equation}{0}

It is under the above understanding of generating Carath\'eodory's definition that we ask which transformations generate the definition. Clearly, for the concept of generation to make sense, any transformation capable of comparing any two sets in $X$ should generate the definition. This is also, in fact, the case, as we show in our main theorem, Theorem $\ref{main}$, that all ergodic $\mu$-outer measure preserving transformations generate the Carath\'eodory definition.

We immediately simplify the task of showing that a transformation generates Carath\'eodory's definition by reducing the 
conditions to be satisfied to properties easier to demonstrate.
\begin{theorem}\label{generatereduction}\thst
Suppose that $(X, \mathcal{B},\mu,T)$ is a dynamical system on a measurable space and let $\mathcal{J}$ be a measure basis. Suppose that for any $T$-invariant $B \subset X$ and for any $J \in \mathcal{J}$, 
$\mu_{\ast}(B \cap K)\geq \mu_{\ast}(B \cap J)$ for each $K \in \mathcal{J}$ with $\mu(K)=\mu(J)$. Then $T$ 
generates the Carath\'eodory definition relative to $\mathcal{J}$.
\pf
Let $B$ be a $T$-invariant set for which there is a $\theta = \theta(B)<2$ for which 
$\mu_{\ast}(B~\cap~J)+\mu_{\ast}(B^c \cap J) < \theta\mu(J)$ for each $J \in \mathcal{J}$. Suppose that 
$\mu_{\ast}(B),\mu_{\ast}(B^c)>0$ and set $\varepsilon=(2-\theta)/2$. Since $\mathcal{J}$ is a measure basis there exist $J_1,J_2 \in \mathcal{J}$ such that 
$$\mu_{\ast}(B \cap J_1)>(1-\varepsilon)\mu(J_1) \hbox{ and }\mu_{\ast}(B^c\cap J_2)>(1-\varepsilon)\mu(J_2),$$ 
for which, by hypothesis,
$$\mu_{\ast}(B \cap J_2)\geq \mu_{\ast}(B \cap J_1) > (1-\varepsilon)\mu(J_1) = (1-\varepsilon)\mu(J_2).$$
That $\mu_{\ast}(B \cap J_2)+\mu_{\ast}(B^c\cap J_2)>2(1-\varepsilon)\mu(J_2) = \theta\mu(J_2)$ follows. This 
contradiction shows that either $\mu_{\ast}(B)$ or $\mu_{\ast}(B^c)$ equals zero, and thus by definition that 
$T$ generates Carath\'eodory's definition.
\end{proof}
\end{theorem}\noindent
\begin{remark} Clearly for any $K \in \mathcal{J}$ we can use $K$ as our original set to get 
$\mu_{\ast}(B \cap J)\geq \mu_{\ast}(B \cap K)$ and thus have $\mu_{\ast}(B \cap K)= \mu_{\ast}(B \cap J)$. However, 
while using `$=$' or `$\leq$' has little effect on the proof of this theorem, the phrasing allows for a simpler proof 
of the main theorem.
\end{remark}

We show the hypothesis of Theorem $\ref{generatereduction}$ and thus the main Theorem, Theorem $\ref{main}$, by showing that  $B\cap K$ can be controllably transported to $B\cap J$.

In the irrational rotation case and the case of Kakutani's transformation, it was possible to track the elements of the measure basis, whose structure and connectedness are sufficiently preserved under the action of the transformations. This is of course not true of a general dynamical system. The difficulty in our main theorem, Theorem $\ref{main}$, is allowing for the first element chosen out of the measure basis to break up under iterations of the transformation and arrive in the second over a range of differing `times' (differing iterations of the transformation). Before proving Theorem $\ref{main}$ we require a simple lemma, whose proof can be found, for e.g., in \cite{Ko1} or \cite{Ni2}.

\begin{lemma}\label{disjointnonmeas}\thst
Suppose that $(X,\mathcal{B},\mu)$ is a measure space, $B \subset X$ and $\{A_n\}_{n\in \N}$ is a sequence of disjoint sets in 
$\mathcal{B}$ such that $\lim_{n \rightarrow \infty}\mu(\cup_{i=n}^{\infty}A_i)=0.$ Then
$$\mu_{\ast}\left(\bigcup_{i=1}^{\infty}(A_i\cap B)\right)=\sum_{i=1}^{\infty}\mu_{\ast}(A_i\cap B).$$
\end{lemma}\noindent
\begin{theorem}\label{main}\thst
Suppose $T$ is an outer measure preserving ergodic transformation on the measure space $(X,\mathcal{B},\mu)$. Then, 
for any measure basis, $\mathcal{J}$, $T$ generates the Carath\'eodory definition relative to $\mathcal{J}$.
\pf
We show that the hypotheses of Theorem $\ref{generatereduction}$ hold. Take 
$J_1, J_2\in \mathcal{J}$ with $\mu(J_1)=\mu(J_2)$. We construct {\it splinters} of $J_1$ with respect to $J_2$ as follows. Let 
$$A_1:=T^{-1}(J_1) \cap J_2 \in \mathcal{B} \ \ \hbox{ and}$$
$$B_1:=T^{-1}(J_1) \sim A_1=T^{-1}(J_1)\sim J_2 \in \mathcal{B},$$
(where here, and in the remainder of the work, $\sim$ denotes set subtraction). We then continue to consider the orbit of $J_1$ allowing parts of elements of the orbit of $J_1$ to eventually cover $J_2$. That is, we define inductively 
$$A_n: = T^{-1}(B_{n-1}) \cap \left(J_2\sim \bigcup_{i=1}^{n-1}A_i\right) \in \mathcal{B}$$
to be the $n$th splinter, 
$$B_n:=T^{-1}(B_{n-1})\sim A_n,$$
and $A_0:= \cup_{i=1}^{\infty}A_i$. We note that the sequence of sets $\{A_i\}_{i=1}^{\infty}$ is necessarily measurable, pairwise disjoint and satisfies
$A_{n+1} \cap \cup_{i=1}^nA_i=\emptyset$
for each $n \in \mathbb{N}$.

We also claim that for each $n \in \mathbb{N}$, $\mu(B_n)=\mu(J_2-\cup_{i=1}^nA_i)$. To prove the claim we observe that
\begin{eqnarray}
\mu(J_2\sim A_1)&=&\mu(J_2)-\mu(A_1)=\mu(T^{-1}(J_1))-\mu(T^{-1}(J_1)\cap J_2) \nonumber \\
&=&\mu(T^{-1}(J_1)\sim J_2)=\mu(B_1) \nonumber
\end{eqnarray}
and, assuming that the claim is true for $n \in \mathbb{N}$, that
\begin{eqnarray}
\mu\left(J_2\sim\bigcup_{i=1}^{n+1}A_i\right)&=&\mu\left(J_2\sim\bigcup_{i=1}^nA_i\right)-\mu(A_{n+1}) \nonumber \\
& = & \mu(B_n)-\mu(A_{n+1}) \nonumber \\
& = & \mu(T^{-1}(B_n)\sim A_{n+1})\nonumber \\
& = & \mu(B_{n+1}), \nonumber
\end{eqnarray}
which shows, through induction, that the claim is true. 
Since $\mu(J_2)<\infty$ and $\{J_2\sim\cup_{i=1}^nA_i\}_i$ is a decreasing sequence of sets it follows that 
$$\lim_{n\rightarrow\infty}\mu(B_n)=\mu\left(J_2\sim A_0\right) \ \ \hbox{ exists.}$$
We now show that $\mu(J_2\sim A_0)=0$. Using that
$B_n\cap J_2\subset J_2 \sim\cup_{i=1}^{n-1}A_i\subset J_2 \sim A_0$ for each $n \in \mathbb{N}$, we first note that for each 
$m,n \in \mathbb{N}$
\begin{eqnarray}
\mu (T^{-m}B_n \cap  (J_2\sim A_0 )) 
& = & \mu \left(\left(B_{n+m}\cup\bigcup_{i=n+1}^{n+m}T^{i-m-n}(A_i)\right) \cap\left(J_2\sim A_0\right)\right) \nonumber \\
& \leq &\mu(B_{n+m}\cap (J_2\sim A_0)) \nonumber \\
& & + \sum_{i=n+1}^{n+m}\mu (T^{i-m-n}(A_i)\cap(J_2\sim A_0)) \nonumber \\
&\leq & \sum_{i=n+1}^{n+m}\mu(T^{i-m-n}A_i) \nonumber \\
& \leq & \mu\left(\bigcup_{i=n+1}^{\infty}A_i\right). \nonumber
\end{eqnarray}
Now, suppose $\lim_{n\rightarrow\infty}\mu(B_n)=\mu(J_2\sim A_0)=:c>0$. Then there exists 
$n_0 \in \mathbb{N}$ such that for each $n \geq n_0$
$$c\leq \mu(B_{n})=\mu\left(J_2\sim\bigcup_{i=1}^nA_i\right)<c+\frac{c^2}{2},$$
and hence
\begin{eqnarray}
\mu (T^{-m}B_{n_0} \cap  (J_2\sim A_0 ) ) & \leq & \mu\left(\bigcup_{i=n_0+1}^{\infty}A_i\right) \nonumber \\
& = & \mu\left(J_2\sim\bigcup_{i=1}^{n_0}A_i\right)-\mu (J_2\sim A_0) \nonumber \\
& < & \frac{c^2}{2} \nonumber
\end{eqnarray}
for each $m \in \mathbb{N}$. Thus
\begin{eqnarray}
\limsup_{m\rightarrow\infty}\frac{1}{m}\sum_{j=1}^m\mu (T^{-j}(B_{n_0})\cap (J_2\sim A_0) )  <  \frac{c^2}{2} < \mu(B_{n_0})\mu (J_2\sim A_0 ), \nonumber
\end{eqnarray}
contradicting $T$'s ergodicity. Thus, $\lim_{n\rightarrow\infty}\mu(B_n)=\mu(J_2\sim A_0)=0.$

Let $B$ be any $T$-invariant set. Then using the outer measure preserving property of $T$ 
\begin{eqnarray}
\mu_{\ast}(J_1\cap B)& = &\mu_{\ast}(T^{-1}(J_1\cap B)) \nonumber \\
& \leq & \mu_{\ast}(J_2\cap T^{-1}(J_1 \cap B))+\mu_{\ast}(T^{-1}(J_1\cap B)\sim J_2) \nonumber \\
& = &\mu_{\ast}((J_2 \cap T^{-1}(J_1))\cap T^{-1}(B))+\mu_{\ast}(B_1\cap B) \nonumber \\
& = & \mu_{\ast}(A_1\cap B)+\mu_{\ast}(B_1 \cap B). \nonumber 
\end{eqnarray}
We continue inductively in a similar manner to see that for each $n \in \mathbb{N}$
$$\mu_{\ast}(J_1\cap B) \leq \mu_{\ast}(B_n\cap B)+\sum_{i=1}^n\mu_{\ast}(A_i \cap B).$$
Also, using the measurability of $A_0$ and Lemma $\ref{disjointnonmeas}$
$$\mu_{\ast}(J_2\cap B)=\mu_{\ast}(A_0\cap B)+\mu_{\ast}((J_2\sim A_0)\cap B) 
=\sum_{i=1}^{\infty}\mu_{\ast}(A_i \cap B).$$
Since $\lim_{n \rightarrow \infty}\mu(B_n)=0$, we choose for any $\varepsilon >0$ an $n_1 \in \mathbb{N}$ such that 
$\mu(B_{n_1})<\varepsilon$ and calculate
$$\mu_{\ast}(J_1\cap B)\leq\mu_{\ast}(B_{n_1}\cap B)+\sum_{i=1}^{n_1}\mu_{\ast}(A_i\cap B)<\varepsilon +\sum_{i=1}^{\infty}
\mu_{\ast}(A_i\cap B) = \varepsilon + \mu_{\ast}(J_2\cap B).$$
Since this is true for each $\varepsilon>0$, $\mu_{\ast}(J_1\cap B) \leq \mu_{\ast}(J_2\cap B)$ fulfilling the 
conditions of Theorem $\ref{generatereduction}$ and thus completing the proof.
\end{proof}
\end{theorem}\noindent
In the next result we refer to a non-trivial measure space by which we mean a measure space, $(X, \mathcal{B},\mu)$, 
in which there exists at least one set $B \in \mathcal{B}$ for which $0<\mu(B)<\mu(X)$
\begin{corollary}\label{mainiff}\thst
Let $T$ be an outer measure preserving transformation on a 
non-trivial measure space $(X,\mathcal{B},\mu)$. Then for any measure basis $\mathcal{J}$, $T$ generates the Carath\'eodory definition relative to $\mathcal{J}$ if and only if $T$ is ergodic.
\pf
That the ergodicity of $T$ ensures that $T$ generates Carath\'eodory's definition follows from Theorem $\ref{main}$. Now 
suppose that $T$ is not ergodic. Then there is a $T$-invariant set $B \in \mathcal{B}$ such that $0<\mu(B)<\mu(X)$. 
As $B$ is 
measurable $\mu(B\cap J)+\mu(B^c\cap J)\leq\theta\mu(J)$ with $\theta=1$ for each $J\in \mathcal{J}$. 
However, $\mu(B) \not=0$ and $\mu(B^c)=\mu(X)-\mu(B)\not=0$. It follows that $T$ does not generate Carath\'eodory's 
definition relative to $\mathcal{J}$.
\end{proof}
\end{corollary}\noindent
Similar results hold for families of transformations with ergodic countable subfamilies of transformations with the 
sacrifice that the transformations (at least of the ergodic subfamily) must be invertible, which is not the case in 
Theorem $\ref{main}$. For details see \cite{Ko1}. 

We consider briefly some properties of the unusual requirement that a transformation be 
outer measure preserving.
\section{Outer measure preserving transformations}
\setcounter{equation}{0}
The properties of outer measure preserving transformations shown below are chosen on the basis of being sufficient to show 
that irrational rotations and Kakutani's transformation are outer measure preserving. We show first that measure 
preserving transformations with measurable images of measurable sets are outer measure preserving. This transfers 
the interest to images of measurable sets. In this direction we make an adjustment to a theorem of Federer \cite{Fe}
from which it is possible to show the appropriate properties of the irrational rotations and Kakutani's transformation.
\begin{theorem}\label{measurableimage}\thst
If $T$ is a measure preserving transformation on $(X,\mathcal{B},\mu)$ such that $T(A)\in \mathcal{B}$ for each 
$A \in \mathcal{B}$, then $T$ is outer measure preserving.
\pf
Let $B\subset X$. Then 
\begin{eqnarray}
\mu_{\ast}(B)&=&\inf\{\mu(A):B\subset A, A\in \mathcal{B}\} \nonumber \\
& = & \inf\{\mu(A):T^{-1}(B)\subset T^{-1}(A),A \in \mathcal{B}\}\nonumber \\
& = & \inf\{\mu(T^{-1}(A)):T^{-1}(B)\subset T^{-1}(A), A\in \mathcal{B}\} \nonumber \\
& = & \inf\{\mu(T^{-1}(A)):T^{-1}(B) \subset T^{-1}(A), T^{-1}(A) \in \mathcal{B}\} \nonumber \\
& = & \mu_{\ast}(T^{-1}(B)). \nonumber
\end{eqnarray}
\end{proof}
\end{theorem}\noindent
The result we adjust to show that the images of measurable sets are measurable under the specific transformations we consider can be stated as in the following theorem. A proof can be found in Federer \cite{Fe}.
\begin{theorem}\label{federer}\thst
Let $(X,\mathcal{B},\nu)$ and $(Y,\mathcal{D},\mu)$ be two measure spaces. Let $X$ be a complete separable metric space 
and $Y$ be a Hausdorff space in which every closed set is $\mu$-measurable. Then, if $f:X \rightarrow Y$ is continuous the 
$f$-image of Borel sets is $\mu$-measurable.
\end{theorem}\noindent
The necessary adjustment to Theorem $\ref{federer}$ required for our results is given below.
\begin{corollary}\label{measimdiscont}\thst
Let $(X,\mathcal{B},\nu)$ and $(Y,\mathcal{D},\mu)$ be two measure spaces. Let $X$ be a complete separable metric space 
and $Y$ be a Hausdorff space in which every closed set is $\mu$-measurable. Then, if $f:X \rightarrow Y$ 
has a discrete set of discontinuities, the $f$-image of Borel sets is $\mu$-measurable.
\pf
Let $D$ be the set of discontinuities. As $D$ is discrete, there exists a 
countable collection of disjoint balls $\{B_d=B_{r_d}(d):d\in D\}$ where $r_d>0$ for each $d\in D$. For any ball 
$B_d$ in this collection let $C_{d,n}=B_d\sim B_{r_d/n}(d)$.

For each $n\in \mathbb{N}$ and $d \in D$ the restriction of $(X,\mathcal{B},\nu)$ to $C_{d,n}$ is a measure space on a 
complete separable metric space and $f$ is continuous on $C_{d,n}$ so that, by Theorem $\ref{federer}$, $f(B \cap C_{d,n})$ 
is $\mu$-measurable for each $B \in \mathcal{B}$. Additionally, since, for each $d\in D$, 
$f(d)$ is a point in $Y$ it is closed and thus $\mu$-measurable.

Moreover, for $X_1=X\sim\cup_{d\in D}B_{r_d}(d)$, the restriction of $(X,\mathcal{B},\nu)$ to $X_1$ is a measure space 
on a complete separable metric space on which $f$ is continuous so that, by Theorem $\ref{federer}$, $f(B\cap X_1)$ is 
$\mu$-measurable for each $B\in \mathcal{B}$. Thus for each $B \in \mathcal{B}$ 
\begin{eqnarray}
f(B) & = & f(B\cap X_1)\cup\bigcup_{d\in D}\left(f\left(B\cap \bigcup_{n\in\mathbb{N}}C_{d,n}\right)\cup f(B\cap \{d\})\right) \nonumber \\
& = & f(B\cap X_1)\cup \bigcup_{d\in D}f(B\cap\{d\})\cup\bigcup_{d\in D}\bigcup_{n\in\mathbb{N}}f(B\cap C_{d,n}), \nonumber
\end{eqnarray}
which is a countable union of $\mu$-measurable sets and thus $\mu$-measurable.
\end{proof}
\end{corollary}\noindent
Finally we show that both irrational rotations and Kakutani's transformation 
are ergodic outer measure preserving transformations. We first give a precise definition of Kakutani's transformation.

Let $X:=[0,1)$ be equipped with the Lebesgue measure, $\mu$, and the usual Borel sets $\mathcal{B}$ to form the measure space $(X,\mathcal{B},\mu)$. We define the {\it primitive} transformation on $X$, $\psi:X\rightarrow X$ by 
$$\psi(x)=x-1+\frac{1}{2^n}+\frac{1}{2^{n+1}} \ \ \hfill x\in I_n$$
where $I_n:= [1-2^{-n},1-2^{n+1})$. Define $A:=\cup_{n=0}^{\infty}I_{2n}$ and let $A^{\prime}$ be a set disjoint from $X$ for which there exists a bijection $\tau:A\rightarrow A^{\prime}$. $A^{\prime}$ is then the first (and in this case only) floor of a so called tower construction (see \cite{Ko1} or \cite{Wa} for further explanation). We then extend the primitive dynamical system $(X,\mathcal{B},\mu ,T)$ as follows: We define the measure space $(\tilde{X},\tilde{\mathcal{B}},\tilde{\mu})$ by 
$$\tilde{X}:=X\cup A^{\prime},$$
$$\tilde{\mathcal{B}}:=\{B:B\subset \tilde{X}, B\cap X \in \mathcal{B} \hbox{ and }\tau^{-1}(B \cap A^{\prime})\in \mathcal{B}\},\hbox{ and}$$
$$\tilde{\mu}(B):=\mu(B\cap X) + \mu(\tau^{-1}(B \cap A^{\prime})),\hbox{ for all }B \in \tilde{B}.$$
We then define Kakutani's transformation $\tilde{\psi}$ on $(\tilde{X},\tilde{\mathcal{B}},\tilde{\mu})$ by
\begin{equation}\label{tower}
 {\widetilde \psi}( {x})=
\begin{cases}
\tau ( {x}),&if \ \ x\in 
A,\cr \psi( {x}),&if \ \ x\in A^c\cap X, \hbox{ and } \cr 
\psi(\tau^{-1}( {x})), &if \ \ x\in A^{\prime},\cr
\end{cases}
\end{equation}

\begin{theorem}\label{irrandkakincluded}\thst
Irrational rotations on the unit circle $\Tbb$ and Kakutani's $\tpsi$ transformation described above are 
ergodic outer measure preserving transformations.
\pf
Since $\rho:\Tbb \rightarrow \Tbb$ is continuous it follows directly from Theorem $\ref{federer}$ that $\rho$ has 
measurable images of measurable sets and thus, by Theorem $\ref{measurableimage}$, is outer measure preserving. As 
mentioned previously, it was shown in \cite{Ni1} that $\rho$ is ergodic.

It has been shown in, for example \cite{Ka1},\cite{Ka2}, and \cite{Ko1} that $\psi$ is ergodic and thus measure preserving. 
Furthermore, we see that the set of points $D=\{1-1/2^n:n\in \mathbb{Z},n\geq 0\}$ is discrete in the usual measure space on $[0,1)$. It follows from the definition of $\psi$ that $D$ is the set 
of discontinuities of $\psi$ and thus, from Corollary $\ref{measimdiscont}$, that $\psi$ is outer measure preserving.

We now show that if a primitive transformation is outer measure preserving, then so to is any 1 level 
tower extension (indeed, the same is true for any tower extension, see \cite{Ko1}, but is not here necessary to show).

Let $B\subset \tilde{X}$. By the definition of a tower extension, ($\ref{tower}$), we see that $\tpsi^{-1}(B\cap X) \subset \tilde{X}\sim A$ and $\tpsi^{-1}(B\cap A^{\prime})\subset A$ so that by the measurability of $A$ 
\begin{eqnarray}
\tilde{\mu_{\ast}(B)}(\tpsi^{-1})& =& \tilde{\mu_{\ast}}(\tpsi^{-1}(B)\cap (\tilde{X}\sim A))+\tilde{\mu_{\ast}}(\tpsi^{-1}(B)\cap A) \nonumber \\
& = & \tilde{\mu_{\ast}}(\tpsi^{-1}(B \cap X))+\tilde{\mu_{\ast}}(\tpsi^{-1}(B \cap A^{\prime}))\nonumber
\end{eqnarray}
and
$$\tilde{\mu_{\ast}}(\tpsi^{-1}(B\cap X))=\tilde{\mu_{\ast}}(\tau(\psi^{-1}(B\cap X)\cap A))+\mu_{\ast}(\psi^{-1}(B\cap X)\sim A).$$
By definition
$\tilde{\mu_{\ast}}(\tau(\psi^{-1}(B\cap X)\cap A))=\mu_{\ast}(\psi^{-1}(B\cap X)\cap A))$ so that, since $\psi$ is outer measure preserving and $A$ is measurable,
\begin{eqnarray}
\tilde{\mu_{\ast}}(\tpsi^{-1}(B\cap X))& = &\mu_{\ast}(\psi^{-1}(B\cap X)\cap A))+\mu_{\ast}(\psi^{-1}(B\cap X)\sim A) \nonumber \\
& = & \mu_{\ast}(\psi^{-1}(B\cap X))=\mu_{\ast}(B\cap X).\nonumber
\end{eqnarray}
Also by definition $\tilde{\mu_{\ast}}(\tpsi^{-1}(B\cap A^{\prime}))=\mu_{\ast}(\tau^{-1}
(B\cap A^{\prime}))=\mu_{\ast}(B\cap A^{\prime})$. Thus
$$\mu_{\ast}(\tpsi^{-1}(B))=\mu_{\ast}(B\cap X)+\mu_{\ast}(B\cap A^{\prime})=\mu_{\ast}(B),$$
and therefore $\tpsi$ is outer measure preserving. The proof is completed by noting that it is known that $\tpsi$ is 
ergodic (see \cite{Ka1}, \cite{Ka2} or \cite{Ko1}).
\end{proof}
\end{theorem} \noindent
\begin{remark} Note that the ergodic (and thus measure preserving) 
transformation on the unit circle, $z\mapsto z^2$ is a continuous function. It follows that the images of measurable sets under this transformation are measurable, and 
thus that the transformation is outer measure preserving. Theorem $\ref{main}$ then shows that this transformation 
does indeed generate the Carath\'eodory definition. The difference here from the analysis in \cite{Ni1} is that, in \cite{Ni1}, the structure of an arc needed to be preserved under negative iterations of the transformation, which is not the case for $z\mapsto z^2$, whereas in the present analysis, by allowing splinters to develop, such structure preservation is not necessary.
\end{remark}
 
\vfill\eject 

\end{document}